\documentclass[12pt, 14paper,reqno]{amsart}
\usepackage{xcolor}
\usepackage{tabularx} 
\usepackage{multirow} 

\usepackage{amsmath,amsfonts,amssymb}
\usepackage[breaklinks]{hyperref}
\usepackage{graphicx}
\usepackage{longtable}
\makeatletter
\@namedef{subjclassname@2020}{%
  \textup{2020} Mathematics Subject Classification}
\makeatother

\theoremstyle{plain}
\newtheorem{thm}{Theorem}[section]
\newtheorem{lem}{Lemma}[section]

\newtheorem{conj}{Conjecture}[section]

\theoremstyle{proof}

\numberwithin{equation}{section}

\begin{document} 
\title[$D(n)$-quadruples in $\mathbb{Z}(\sqrt{4k+2})$]{On a conjecture of Franu\v si\'c and Jadrijevi\' c: Counter-examples}
\author{Kalyan Chakraborty, Shubham Gupta and Azizul Hoque}
\address{K. Chakraborty @Kerala School of Mathematics, Kozhikode-673571, Kerala, India.}
\email{kalychak@ksom.res.in}
\address{S. Gupta @Harish-Chandra Research Institute,  A CI of Homi Bhabha National
Institute, Chhatnag Road, Jhunsi, Prayagraj - 211019.}
\email{shubhamgupta@hri.res.in}
\address{A. Hoque @Department of Mathematics, Rangapara College, Rangapara, Sonitpur-784505, Assam, India.}
\email{ahoque.ms@gmail.com}
\date{\today}
\maketitle

\begin{abstract} Let $d\equiv 2\pmod 4$ be a square-free integer such that $x^2 - dy^2 =- 1$ and $x^2 - dy^2 =  6$ are solvable in integers. We prove the existence of  infinitely many quadruples in $\mathbb{Z}[\sqrt{d}]$ with the property $D(n)$ when $n \in \{(4m + 1) + 4k\sqrt{d}, (4m + 1) + (4k + 2)\sqrt{d}, (4m + 3) + 4k\sqrt{d}, (4m + 3) + (4k + 2)\sqrt{d}, (4m + 2) + (4k + 2)\sqrt{d}\}$ for $m, k \in \mathbb{Z}$.
As a consequence, we provide few counter examples to a conjecture of Franu\v si\'c and Jadrijevi\' c (see Conjecture \ref{Con1.1}).
\end{abstract}

\section{Introduction}
Given a non-zero integer $n$, a set of $m$ distinct positive integers $\{a_1, a_2, \cdots, a_m\}$ such that $a_ia_j + n$ is a perfect square for all $1\leq i<j\leq m$, is called a Diophantine $m$-tuple with the property $D(n)$ or simply $D(n)$-$m$-set. The $D(1)$-$m$-sets are known as Diophantine $m$-tuples, and have been studied since the time of Diophantus.
The first Diophantine quadruple, viz. $\{1, 3, 8, 120\}$ was found by Fermat.
It follows from the work of Baker and Davenport \cite{BD1969} that the above quadruple cannot be extended to a Diophantine quintuple. In 2004, Dujella \cite{DU2004} proved that there are only finitely many Diophantine quintuples. However, there was a `{\it folklore}' conjecture which states that there does not exist Diophantine quintuples. Recently this conjecture has been settled by He, Togb\'e and Ziegler in \cite{HTZ2019}. On the other hand, Bonciocat et al. \cite{BCM2020} settled the conjecture of Dujella \cite{DU1993} which states that there is no Diophantine $D(-1)$-quadruple. We direct the readers to \cite{BR1985, CGH-2022, DU21, EFF2014, HT2011} for more information about $D(n)$-$m$-sets.  

A general notion is obtained by considering elements of any commutative ring instead of positive integers. However, many demanding and interesting problems already occur in the setting fixed by the above definition. For a given non-zero element $n$ of a commutative ring $\mathcal{R}$ with unity, a set $\{a_1, a_2,\cdots, a_m\} \subset \mathcal{R}\setminus\{0\}$ is called a Diophantine $m$-tuple with the property $D(n)$ or simply $D(n)$-$m$-set in $\mathcal{R}$ if 
$
a_i a_j+n
$
is a perfect square in $\mathcal{R}$ for  all $1 \leq i < j \leq m$. 

In \cite{DU1993, MR1985}, it was established that there exists a quadruple with the property $D(n)$ if and only if $n$ can be expressed as the difference of two squares, up to finitely many exceptions. Later, Dujella \cite{DU1997} proved that this fact also true in case of Gaussian integers. Besides this, it was proved that for certain integers $d$, the above fact is true in the ring of integers of  $\mathbb{Q}(\sqrt{d})$ for certain $d$'s (see, \cite{MR2004, FR2004, FR2008, FR2009, FS2014, SO2013}). These results motivate Franu\v{s}i\'c and Jadrijevi\'c to post the following conjecture.

\begin{conj}[{\cite[Conjecture 1]{FJ2019}}]\label{Con1.1}
Let $\mathcal{R}$ be a commutative ring with unity $1$ and $n \in \mathcal{R}\setminus \{0\}$. Then there exists a quadruple with the property $D(n)$ if and only if  $n$ can be expressed as a difference of two squares in $\mathcal{R}$, up to finitely many exceptions.
\end{conj}
This conjecture was verified successfully for ring of integers of certain quadratic fields (see, \cite{MR2004, FR2004, FR2008, FR2009, FS2014, SO2013}). It follows from \cite{FR2013} and \cite{MA2012} that Conjecture \ref{Con1.1} is true for the ring of integers of $\mathbb{Q}(\sqrt[3]{2})$. Recently Franu\v{s}i\'c and Jadrijevi\'c \cite{FJ2019} confirmed this conjecture for the ring of integers of bicyclic biquadratic field $\mathbb{Q}[\sqrt{2}, \sqrt{3}]$.  

In this paper,  we investigate the existence of infinitely many quadruples in $\mathbb{Z}[\sqrt{d}]$ with the property $D(n)$ for certain $n\in\mathbb{Z}[\sqrt{d}]$ under the assumptions 
that the equations 
\begin{equation}\label{eqi1}
x^2 - dy^2 = -  1
\end{equation}
and 
\begin{equation}\label{eqi2}
x^2 - dy^2 =   6
\end{equation}
are solvable in integers. Precisely, we prove the following:

\begin{thm}\label{1.1}
Assume that $d\equiv 2\pmod 4$ is a square-free positive integer and the equations \eqref{eqi1} and \eqref{eqi2} are solvable. Then there exist infinity many quadruples in $\mathbb{Z}[\sqrt{d}]$ with the property $D(n)$ when $n \in \{(4m + 1) + 4k\sqrt{d}, (4m + 1) + (4k + 2)\sqrt{d}, (4m + 3) + 4k\sqrt{d}, (4m + 3) + (4k + 2)\sqrt{d}, (4m + 2) + (4k + 2)\sqrt{d}\}$ with $m, k \in \mathbb{Z}$. 
\end{thm}

As a consequence of Theorem \ref{1.1}, we construct some counter-examples of Conjecture \ref{Con1.1}.  In fact, under the assumption of Bunyakovsky's conjecture, we see that there are infinitely many such rings in which Conjecture \ref{Con1.1} does not hold.
\section{Preliminary Descent}
We first fix up the following notations:
\begin{itemize}
\item $
(a, b) := a + b\sqrt{d}$,
\item $
Nm(\alpha) := (a, b)(a, -b)$ for a given $\alpha=(a, b)$,
\item $(x, y) \equiv (a, b) \text{~(mod~} (c, d))$ denotes
$
x \equiv a \pmod{c} \text{~~~and~~~} y \equiv b \pmod{d}.
$
\end{itemize}
For any square-free integer $d \equiv 2 \pmod{4}$ and $m, k \in \mathbb{Z}$, we define:
\begin{align*}
\mathcal{S} :=& \{(4m, 4k + 1), (4m, 4k + 2), (4m, 4k + 3), (4m + 1, 4k + 1), (4m + 1, 4k + 3), (4m + 2, \\ & 4k + 1), (4m + 2, 4k + 3), (4m + 3, 4k + 1),(4m + 3, 4k + 3)\},\\ 
\mathcal{T} :=& \{(4m, 4k), (4m + 1, 4k), (4m + 1, 4k + 2), (4m + 2, 4k), (4m + 2, 4k + 2), (4m + 3, 4k),\\ & (4m + 3, 4k + 2)\}.
\end{align*}
Clearly, if $n \in \mathbb{Z}[\sqrt{d}]$ then $n \in \mathcal{S} \cup \mathcal{T}$. We now extract the following lemma from \cite[Proposition 5]{FR2004}.

\begin{lem}\label{thm2.1} 
Let $d\in \mathbb{Z}$ such that $d\equiv 2\pmod 4$. Then for $n \in \mathcal{S}$, there does not exist any Diophantine quadruple in $\mathbb{Z}[\sqrt{d}]$ with the property $D(n)$.
\end{lem}
There are many results concerning the integers which can be represented as $n=\alpha^2-\beta^2 $ in $\mathbb{Z}[\sqrt{d}]$. These integers 
can be used in the construction of quadruples in $\mathbb{Z}[\sqrt{d}]$. In this context, we deduce the following lemma from \cite[Propositions 5 and 6]{DF2007}.
\begin{lem}\label{thm2.2} 
For any square-free $d\in \mathbb{Z}$, all the numbers of the form $n = (4m + 2, 4n)\in \mathbb{Z}[\sqrt{d}]$ can be written as the difference of two squares if and only if $x^2- dy^2 = \pm 2$ is solvable in 
integers. Further, all the numbers of the form $n = (4m + 2, 4n + 2)\in\mathbb{Z}[\sqrt{d}]$ can be written as the difference of two squares  if and only if $d \equiv 2 \pmod{4}$ and  $x^2 - dy^2 = \pm 2$ is solvable in integers.
\end{lem}

The next lemma follows from the definition. 
\begin{lem}\label{lemp}
For a square-free $d\in \mathbb{Z}$, let  $\{a_1, a_2, a_3, a_4\}$ be a quadruple with the property $D(n)$. Then for any non-zero $w\in \mathbb{Z}[\sqrt{d}]$,  $\{wa_1, wa_2, wa_3, wa_4\}$ is a quadruple in $\mathbb{Z}[\sqrt{d}]$ with the property $D(w^2n)$.
\end{lem}
Let us consider the set $\{a, b, a + b + 2r, a + 4b + 4r\}$ having non-zero distinct elements of $\mathbb{Z}[\sqrt{d}]$. The condition `non-zero distinctness' in the above set can be verified using the following lemma.
\begin{lem}\label{lemm1} Let $a_1, b_1, c_1, d_1, e_1, a_2, b_2, c_2, d_2 \in \mathbb{Z}$ with $a_1, b_1, a_2 \neq 0$. Then
the following system of equations
\begin{align}\label{eQ7}
\begin{cases} a_1x^2 + b_1y^2 + c_1x + d_1y + e_1 = 0, \\
a_2xy + b_2x + c_2y + d_2 = 0,
\end{cases}
\end{align}
\end{lem}
has only finitely many integer solutions. 

\begin{proof}
Consider $a_2xy + b_2x + c_2y + d_2 = 0$. This implies that 
\begin{equation}\label{eqqq1}
x = -\dfrac{c_2y+ d_2}{a_2y + b_2},
\end{equation}
where $y \neq -\dfrac{b_2}{a_2}$. If $y \neq -\dfrac{b_2}{a_2}$, we put this value of $x$ in the first equation of \eqref{eQ7}. This gives a fourth degree equation in $y$ with the leading coefficient $a_2^2b_1 \neq 0$. Thus, we get at most four values of $y$ and therefore by \eqref{eqqq1}, one gets at most four values of $x$. On the other hand, for $y = -\dfrac{b_2}{a_2}$, the first equation of \eqref{eQ7} gives a two degree equation in $x$  with leading coefficient $a_2^2a_1$ which provides at most two values for $x$. This completes the proof.
\end{proof}

We now assume that $n\in \mathbb{Z}[\sqrt{d}]$ such that $ab + n = r^2 $. Then 
$a(a + b + 2r)+n=(a+r)^2$, $b(a + b + 2r)+n=(b+r)^2$ and $b(a + 4b + 4r)+n=(2b+r)^2$. Now if $a(a + 4b + 4r) + n = \alpha^2$ for some  $\alpha  \in \mathbb{Z}[\sqrt{d}]$, then $a^2+4ar+4r^2= (a + 2r)^2 = \alpha^2+3n$. Thus, we have the following lemma.

\begin{lem}\label{lem2.1}
For any $n\in\mathbb{Z}[\sqrt{d}]$, the set $\{a, b, a + b + 2r, a + 4b + 4r\}$ is a quadruple in $\mathbb{Z}[\sqrt{d}]$ with the property $D(n)$ when 
$ ab + n = r^2$ and $3n = \alpha_1\alpha_2$, where $\alpha_1=a + 2r + \alpha$ and $ \alpha_2= a + 2r - \alpha$
for some $ a, b, r, \alpha\in \mathbb{Z}[\sqrt{d}]$ such that all the element of this set are non-zero and distinct.
\end{lem}

\section{Pell type equations $x^2-dy^2=\pm 6$ and $x^2-dy^2=\pm 1$}\label{S1}
We now consider  only  quadratic ring of integers $\mathbb{Z}[\sqrt{d}]$, where the square-free integer $d$ satisfies $d \equiv 2 \pmod 4$, and the equations $x^2 - dy^2 = - 1$, $x^2 - dy^2 =  6$ are solvable in integers. This implies that the equation $x^2 - dy^2 = -6$ will be also solvable in integers. Also it is well known that the ring of integers of real quadratic field has infinitely many elements of norm 1. Now, we discuss about the elements of norm $\pm 1$ and $ \pm 6$ in $\mathbb{Z}[\sqrt{d}]$. First, let $x^2 - dy^2 = 6$ be  solvable in integers, then
$
x^2 - dy^2 \equiv 0 \pmod{3}.
$ 
Clearly, $x$ and $y$ can not be $0$ modulo $3$, since $9$ can not  
divide $6$. If $x \equiv 1 \pmod 3$ and $y \equiv 0 
\pmod3$, then $1 \equiv 0 \pmod 3$, which is not 
possible, and similarly, $x \equiv 2 \pmod3$ and $y \equiv 0 \pmod3$ is not possible. Now, assume that $x \equiv 0 \pmod3$ and $y \equiv \pm 1 \pmod3$. Then  $3$ divides $d$ and thus by the assumption $d \equiv 2 \pmod 4$ one gets $d \equiv 6 \pmod {12}.$  Since the equation $x^2 - dy^2 = -1$ is solvable in integers, so that reading it modulo $12$ we see that  $x \equiv 0 \pmod3$ and $y \equiv \pm 1 \pmod3$ are not possible.
Therefore $x, y \equiv \pm 1 \pmod 3$ and, after reducing  the equation $x^2 - dy^2 = 6$ to modulo 4, we see that $x$ is even and $y$ is odd. This implies that $(x, y)  = (4, 1), (4, 5), (2, 1), (2, 5) \pmod{(6, 6)}$.
Taking $x, y \equiv \pm 1 \pmod 3$, one gets $d\equiv 1\pmod 3$. Therefore $d \equiv 10 \pmod {12}$ since by hypothesis $d \equiv 2 \pmod 4$. 

Consider the Pell equation $x^2 - dy^2 = 1$. After using $d \equiv 10 \pmod{12}$ and reducing this equation to modulo $3$ and modulo $4$, we get $x \equiv \pm 1 \pmod {6}$ and $y \equiv 0 \pmod{6}$.
Similarly for $x^2 - dy^2 = -1$, we get $x \equiv  3 \pmod {6}$ and $y \equiv \pm 1 \pmod{6}$. We also know that there exist infinitely many solutions of these two equations.

Let $\alpha = \beta\gamma$, where $\beta = (6a + 4, 6b + 1)$ and $\gamma = (6a_1 + 3, 6b_1 + 1)$ for some $a, b, a_1, b_1 \in\mathbb{Z}$, with norm of $\beta$ and $\gamma$ are $\pm 6$ and -1, respectively. If $a$ is even then $\beta = (12a + 4, 6b + 1)$, for some $a, b \in \mathbb{Z}$ with norm $\pm 6$, and $\alpha = (12a_1 - 2, 6b_1 + 1)$ with norm $\mp 6$, for some $a_1, b_1 \in \mathbb{Z}$.  In the similar way,  if 
$a$ is odd, then $\alpha = (12M  + 4, 6N + 1)$ for some $M, N \in \mathbb{Z}$ with norm $\mp 6$. Multiplying this $\alpha$ by $\gamma$ one gets 
an element of the form $(12M + 2, 6N + 1)$ with norm $\pm 6$ for some $M, N\in \mathbb{Z}$ because if $x^2-dy^2 = \pm 6$ has a solution $(X, Y)$ then $(\pm X, \pm Y)$ are also solutions.
Now, we prove that there does not exist any norm $-6$ element of the form $(12M \pm 4, 6N \pm 1)$ for some $M, N \in \mathbb{Z}$. Suppose it exists then the equation $(12M \pm 4)^2 - d(6N \pm 1)^2 = -6$ gives $d \equiv 22 \pmod{24}$. In the above paragraph, we have proved that all elements whose norm are $-1$ have the form $(6a \pm 3, 6b \pm 1)$ with $a, b \in \mathbb{Z}$. So that $(6a \pm 3)^2 - d(6b \pm 1)^2  = -1$. After rearranging this equation we get $(36(a)(a + 1) + 9) - d(24b^2 + 12b(b + 1) + 1) = - 1$. Using $d \equiv 22 \pmod{24}$, we get a contradiction. To conclude, in these quadratic ring of integers, norm $6$ and norm $-6$ element will be  of form $(12M \pm 4, 6N \pm 1)$
 and $(12M \pm 2, 6N \pm 1)$ respectively for some $M, N \in \mathbb{Z}$, i.e. $(12M \pm 4)^2 - d(6N \pm 1)^2 = 6$. And if we reduce this equation modulo $48$, we get $d \equiv 10 \pmod{48}$, since $d \equiv 2 \pmod{4}$.

Now we prove that there are infinitely many elements of norm $\pm 6$. Let $\alpha = \beta\gamma$, where $\beta = (12a + 2, 6b + 1)$ and $\gamma = (6a_1 + 1, 6b_1)$ for some $a, b, a_1, b_1 \in \mathbb{Z}$, with norm of $\beta$ and $\gamma$ are $- 6$ and $1$ respectively. Then $\alpha$ will be of the form $(12M + 2, 6N + 1)$, for some $M, N \in \mathbb{Z}$ with norm $- 6$. Since we have infinitely many choices of $\gamma$, there exist infinitely many elements of norm $- 6$ and of the form $(12M \pm 2, 6N \pm 1)$, for $M, N \in \mathbb{Z}$ because if $x^2 - dy^2 = -6$ has a solution $(X, Y)$, then $(\pm X, \pm Y)$ are also solutions. Similarly, there exist infinitely many elements of norm $ 6$ and of the form $(12M \pm 4, 6N \pm 1)$ for some $M, N \in \mathbb{Z}$. We sum up the above discussion in the following lemma.

\begin{lem}\label{lem3.1} Assume that $d\equiv 2\pmod 4$ is a square-free integer and the equations \eqref{eqi1} and \eqref{eqi2} are solvable in integers. Then   
for the ring of integers $\mathbb{Z}[\sqrt{d}]$, the following statements hold:
\begin{itemize}
\item[(i)] all elements whose norm is $1$ have the form $(6a_1 \pm 1, 6b_1)$ and there are infinitely many such elements;
\item[(ii)] all elements whose norm is $-1$ have the form $(6a \pm 3, 6b \pm 1)$ and there are infinitely many such elements;
\item[(iii)] $d \equiv 10 \pmod{48}$;
\item[(iv)] all elements whose norm is $6$ have the form $(12M \pm 4, 6N \pm 1)$ and there are infinitely many such elements;
\item[(v)] all elements whose norm is $-6$ have the form $(12M \pm 2, 6N \pm 1)$ and there are infinitely many such elements,
where $a, b, a_1, b_1, M, N \in \mathbb{Z}$.
\end{itemize}
\end{lem}

\section{Proof of Theorem \ref{1.1}}
The proof is divided in various parts depending on the values of $n$ and it relies upon Lemmas \ref{lem3.1} and \ref{lem2.1}. Lemma \ref{lem3.1} helps us for choosing elements with specific norms whereas Lemma \ref{lem2.1} helps us to construct quadruples. 
\subsection*{ Case 1: $n = (4m + 1, 4k)$} In this case $3n = (-1)(-3)(4m + 1, 4k)$ and thus by (i) of Lemma \ref{lem3.1}, we have 
\begin{align*}
3n &= (-6\alpha - 3, 6\beta + 3)(2\alpha + 1, 2\beta + 1)(4m + 1, 4k), \text{ for some } \alpha, \beta \in \mathbb{Z}\\
&=(-6\alpha - 3, 6\beta + 3)(8m\alpha + 2\alpha + 4m + 1 + d(8\beta k + 4k), 8\alpha k + 4k + 8\beta m + 2\beta + 4m + 1).
\end{align*}
Let us set:
$$\begin{cases}
\alpha_1 = (-6\alpha - 3, 6\beta + 3)\\
\alpha_2 = (8m\alpha + 2\alpha + 4m + 1 + d(8\beta k + 4k), 8\alpha k + 4k + 8\beta m + 2\beta + 4m + 1).
\end{cases}$$
Then by Lemma \ref{lem2.1}, we have
$$\alpha_1+\alpha_2=2a+4r,$$ which implies that 
\begin{equation}\label{4.1}
a + 2r = (4m\alpha - 2\alpha + 2m - 1 + 2d(2\beta k + k), 4\alpha k + 2k + 4\beta m + 4\beta + 2m + 2).
\end{equation}

Now by  (i) of Lemma \ref{lem3.1}, we have infinitely many $a\in \mathbb{Z}[\sqrt{d}]$ such that $a = (6a_1 + 1, 6b_1)$ with $Nm(a)=1$.  Note that out of these infinitely many choices of $a$, by Lemma \ref{lemm1}, only finitely many provide a zero and distinct element in the set $\{a, b, a+b+2r, a+4b+4r\}$ as in Lemma \ref{lem2.1}. Thus \eqref{4.1} gives 
$$
r = (2m\alpha - \alpha + m - 1 -3a_1 + d(2\beta k + k), 2\alpha k + k + 2\beta m + 2\beta + m + 1 - 3b_1),$$
which is
$$
r=(m(2\alpha + 1) + kd(2\beta + 1) - \alpha -1 - 3a_1, m(2\beta + 1) + k(2\alpha + 1) + 2\beta + 1 - 3b_1).
$$
Again from Lemma \ref{lem2.1}, we use $ab + n = r^2$ to calculate $b$ which turned out to be 
\begin{align*}
b =& \Big((m(2\alpha + 1) + kd(2\beta + 1) - \alpha -1 - 3a_1)^2 + d(m(2\beta + 1) + k(2\alpha + 1) + 2\beta + 1 -\\
&  3b_1)^2 -4m - 1, 2(m(2\alpha + 1) + kd(2\beta + 1) - \alpha -1 - 3a_1)(m(2\beta + 1) + k(2\alpha + 1) + \\ & 2\beta + 1 - 3b_1) - 4k\Big)\Big(6a_1 + 1, -6b_1\Big).
\end{align*}
Thus by Lemma \ref{lem2.1}, the elements $a, b$ and $r$ give the quadruple $\{a, b,a+b+2r,a+4b+4r \}$ in $\mathbb{Z}[\sqrt{d}]$ with the property $D(4m + 1, 4k)$. Since there are infinitely many choices for $a$, we conclude that there are infinitely many quadruples of the form  $\{a, b,a+b+2r,a+4b+4r \}$ in $\mathbb{Z}[\sqrt{d}]$ with the property $D(4m + 1, 4k)$. 

\subsection*{Case 2: $n = (4m + 1, 4k + 2)$} Here, we will use the same $\alpha_1$ and $\alpha_2$ as in the previous case.
Analogous to \eqref{4.1}, we get
$$a + 2r = \Big(4m\alpha - 2\alpha + 2m -1 + d(4\beta k + 2\beta + 2k + 1),
4\alpha k + 2\alpha + 2k + 3 + 4m\beta + 4\beta + 2m\Big).$$
Now we choose $a = (6a_1 + 3, 6b_1 + 1)$ with $Nm(a)=-1$ and the above equation gives 
$$r = \Big(2m\alpha - \alpha + m -2 - 3a_1 + (d/2)(4\beta k + 2\beta + 2k + 1),
2\alpha k + \alpha + k + 1 - 3b_1 + 2m\beta + 2\beta + m\Big).$$

We complete this case by computing $b$ using $ab + n = r^2$, which is given by 
\begin{align*}
b = &\Big((2m\alpha - \alpha + m -2 - 3a_1 + (d/2)(4\beta k + 2\beta + 2k + 1)^2 + d(2\alpha k + \alpha + k + 1 - 3b_1 +\\
& 2m\beta + 2\beta + m)^2 - 4m - 1,  2(2m\alpha - \alpha + m -2 - 3a_1 + (d/2)(4\beta k + 2\beta + 2k + 1)\\
& (2\alpha k + \alpha + k + 1 - 3b_1 + 2m\beta + 2\beta + m) - 4k - 2\Big) \times \Big(-6a_1 - 3, 6b_1 + 1\Big).
\end{align*}

\subsection*{Case 3: $n = (4m + 3, 4k)$}
In this case, $3n = 3(4m + 3, 4k)$ and we consider $\alpha_1 = 3$ and  $\alpha_2 = (4m + 3, 4k)$. Using Lemma \ref{lem2.1}, we get
$$
a + 2r = (2m + 3, 2k).
$$

Employing Lemma \ref{lem3.1}, we can find infinitely many elements in $\mathbb{Z}[\sqrt{d}]$ of the form $$a = (6a_1 + 1, 6b_1)$$ with $Nm(a)=1$, and thus the last equation gives
$$
r = (m + 1 - 3a_1, k - 3b_1).
$$
Finally $ab + n = r^2$ gives 
$$
b =\Big\{(m + 1 - 3a_1)^2 + d(k - 3b_1)^2 - 4m - 3, 2(m + 1 - 3a_1)(k - 3b_1) - 4k\Big\}\Big(6a_1 + 1, -6b_1\Big). 
$$
\subsection*{Case 4: $n = (4m + 3, 4k + 2)$}
Let $\alpha_1$ and $\alpha_2$ be as in Case 3. Then by Lemma \ref{lem2.1}, we have
$
a + 2r = (2m + 3,~2k + 1). 
$
We choose $$a = (6a_1 + 3,~6b_1 + 1)$$ in $\mathbb{Z}[\sqrt{d}]$ with $Nm(a)=-1$. Then
$$
r = (m - 3a_1, k -3b_1).
$$
Again the equation $ab + n =r^2$ gives
$$
b = \Big\{(m - 3a_1)^2 + d(k -3b_1)^2 - 4m - 3,~ 2(m - 3a_1)(k -3b_1) - 4k - 2\Big\}\Big(-6a_1 - 3, 6b_1 + 1\Big).
$$

\subsection*{Case 5: $n = (4m + 2, 4k + 2)$} Here, we have 
\begin{align}\label{eqq5}
3n &= 3(4m + 2, 4k + 2) \nonumber\\ &= 6(2m + 1, 2k + 1) \nonumber\\
& = (12M + 4, -6N - 1)(12M + 4, 6N + 1)(2m + 1, 2k + 1) \text{ for some }M, N\in \mathbb{Z}.
\end{align}

Now we choose $\alpha_1$, $\alpha_2$ as follows:
\begin{align*}
\alpha_ 1 &= (12M + 4, -6N - 1),\\
\alpha_2 &= (12M + 4, 6N + 1)(2m + 1, 2k + 1) \\ 
&= (24Mm + 12M + 8m + 4 + d(12Nk + 6N + 2k + 1), 24Mk + 12M + 8k + 5\\ & +  12Nm + 6N + 2m).
\end{align*}
As before, utilizing Lemma \ref{lem2.1}, we obtain
\begin{equation}\label{eQ1}
a + 2r = (12Mm + 12M + 4m + 4 + \frac{d}{2}(12Nk + 6N + 2k + 1), 12Mk + 6M + 4k + 2 + 6Nm + m).
\end{equation}
Again by Lemma \ref{lem3.1}, we can find infinitely many $a = (6a_1 + 1, 6b_1)$ in $\mathbb{Z}[\sqrt{d}]$ with $Nm(a)=1$, and thus \eqref{eQ1} gives 
\begin{align*}
r = &(6Mm + 6M + 2m + 2 + \frac{d}{2}(6Nk + 3N + k) + \frac{1}{2}(\frac{d}{2} - 1) - 3a_1,
6Mk + 3M +\\
&  2k + 1 + 3Nm + \frac{m}{2} - 3b_1).
\end{align*} 
For $r \in \mathbb{Z}[\sqrt{d}]$, $m$ should be even, and thus the equation $ab + n  = r^2$ entails
\begin{align*}
b = & \Big((6Mm + 6M + 2m + 2 + (d/2)(6Nk + 3N + k) + (1/2)(d/2 - 1) - 3a_1)^2 +  d(6Mk +\\ & 3M + 2k + 1 + 3Nm + (m/2) - 3b_1)^2 - 4m - 2, 
 2(6Mm + 6M +  2m + 2 + (d/2)(6Nk \\ &+ 3N + k) + (1/2)(d/2 - 1) - 3a_1)
(6Mk + 3M + 2k + 1 + 3Nm + (m/2) - 3b_1) - 4k\\ & - 2\Big)\times \Big(6a_1 + 1, -6b_1\Big).
\end{align*}

On the other hand for odd $m$, we choose $a = (6a_1 + 3, 6b_1 + 1)$ with $Nm(a)=-1$. In this case \eqref{eQ1} gives
\begin{align*}
r = & \Big(6Mm + 6M + 2m + 2 + (d/2)(6Nk + 3N + k) + (1/2)(d/2 - 3) - 3a_1,
6Mk + 3M \\  & + 2k + 1 + 3Nm + (m/2) - 3b_1 - (1/2)\Big),
\end{align*}
and hence
\begin{align*}
b = &\Big((6Mm + 6M + 2m + 2 + (d/2)(6Nk + 3N + k) + (1/2)(d/2 - 3) - 3a_1)^2 +  d(6Mk +\\  & 3M + 2k + 1 + 3Nm + (m/2) - 3b_1 - (1/2))^2 - 4m - 2,  2(6Mm + 6M + 2m + 2 + (d/2)\\  &(6Nk + 3N + k) + (1/2)(d/2 - 3) - 3a_1) (6Mk + 3M + 2k + 1 + 3Nm + (m/2) - 3b_1 - \\  & (1/2)) - 4k - 2\Big) \times \Big(-6a_1 - 3, 6b_1 + 1\Big).
\end{align*}
This completes the proof of Theorem \ref{1.1}.

\section{Counter-examples of Conjecture \ref{Con1.1}}
We have $d\equiv 10 \pmod {48}$ by Lemma \ref{lem3.1}. Thus utilizing Lemma \ref{thm2.2} we see that all the numbers of the form $(4m + 2) + (4k + 2)\sqrt{d}$  can be expressed as the difference of two square in $\mathbb{Z}[\sqrt{d}]$ if and only if $x^2 - dy^2 = \pm 2$ is solvable in integers.  We now consider  the equation $x^2 - dy^2 = \pm 2$. Reducing this equation modulo $48$, we see that no integers $x$ and $y$ satisfy it. To conclude, we state the following lemma.
\begin{lem}\label{lem5.1} Let $d \equiv 2 \pmod4$ be square-free integer. Assume that the ring $\mathbb{Z}[\sqrt{d}]$ contains elements of norm $-1$ and norm $6$. Then $\mathbb{Z}[\sqrt{d}]$ does not contain any element of norm $\pm 2$.
\end{lem}
Sum up these, one sees that for $d \equiv 10 \pmod{48}$, there are integers $m$ and $k$ such that $(4m + 2) + (4k + 2)\sqrt{d}$ can not be expressed as the difference of two square in $\mathbb{Z}[\sqrt{d}]$. Further we will prove in fact that there exist infinitely many $n = (4m + 2) + (4k + 2)\sqrt{d}$  which are not representable as a difference of two squares in $\mathbb{Z}[\sqrt{d}]$. Then by employing Theorem \ref{1.1} on these $n$'s, we get that a Diophantine quadruple in $\mathbb{Z}[\sqrt{d}]$ with the property $D(n)$ exists but $n$ is not expressible as the difference of two squares in $\mathbb{Z}[\sqrt{d}]$. 

Assume that $n = (4m + 2) + (4k + 2)\sqrt{d}$.  Dirichlet announced that among the primes represented by the quadratic form $ax^2 + 2bxy + cy^2$, where $\gcd (a, 2b, c) = 1$, infinitely many of them are representable by any given linear form $Mx+N$, with $\gcd(M,N) = 1$, provided $a,b,c,M,N$ are such that the linear and quadratic forms can represent the same number. Meyer \cite{ME1888} gave a complete proof of this result, whereas  Mertens \cite{ME1895} gave an elementary proof of the same. We conclude by this result that  there exist infinitely many  primes $p$ of the form $p \equiv 3 \pmod{4}$ such that 
$
x^2-dy^2=p.
$ 
As $d\equiv 2\pmod 4$, so that both $x$ and $y$ are odd. Thus, there exist infinitely many  primes $p$ of the form $p \equiv 3 \pmod{4}$ such that 
$
(2m + 1)^2 - d(2k + 1)^2 = p.
$ 
Therefore due to the existence of infinitely many such primes, we get infinitely many pairs of integers $m$ and $k$. We claim that for these $m$ and $k$,  $n = (4m + 2) + (4k + 2)\sqrt{d}$ can not be written as difference of two squares. To prove this claim, we follow  \cite[Proposition 6, (ii)]{DF2007}, which gives that $u^2 - dv^2 = \pm 2$ must be solvable in integers $u, v$. But from Lemma \ref{lem5.1}, it is not possible. Hence  there exist infinitely many $n$ which can not be written as a difference of two squares in $\mathbb{Z}[\sqrt{d}]$, however Theorem \ref{1.1} confirms the existence of a Diophantine quadruple in $\mathbb{Z}[\sqrt{d}]$ with the property $D(n)$. 

Constructing a counter-example for Conjecture \ref{Con1.1}, we can take square-free positive integer $d \equiv 2 \pmod{4}$ such that $x^2 - dy^2 = 6$ and $x^2 - dy^2 = -1$ are solvable in integers $x$ and $y$. Let $d = 48l + 10$ with $l \in \mathbb{N} \cup \{0\}$ (see Lemma 3.1). In the following table, we record some positive integer $l$ for which \eqref{eqi1} and \eqref{eqi2} are solvable: 

\begin{center}
\begin{tabular}{ |p{0.5cm}|p{.5cm}|p{4.5cm}|p{4.5cm}| }
\hline
$l$ & $d$ & $(x, y)$ for \eqref{eqi1} & $(x, y)$ for \eqref{eqi2} \\
\hline
0 & 10  & (3,1) & (4,1)\\
1 & 58   & (99,13) & (8,1)\\
2 & 106 & (4005,389) & (1184,115) \\
4 & 202 & (3141,221) & (668,47) \\
6 & 298 & (409557,23725) & (328,19) \\
\hline
\end{tabular}
\end{center}
\vspace{0.5cm}
We now determine $l$ for which $x^2 - dy^2 = 6$ and $x^2 - dy^2 = -1$ are solvable in integers $x$ and $y$. First we consider $x^2 - (48l + 10)y^2 = 6$. Putting $y = 1$ and $l = 3l'^2 \pm 2l'$ in this equation, we get $x = (12l'  \pm 4)$, where $l' \in \mathbb{Z}$. Hence for $d = 48(3l'^2 \pm 2l') + 10 = 2(24l'(3l' \pm 2) + 5)$, the equation $x^2 - dy^2 = 6$ will be solvable. If $d = 2p$ , where $p$ is prime and $p \equiv 5 \pmod{8}$, then $x^2 - 2py^2 = -1$ is solvable in integers (see \cite{LEM}). Let $p(l') = 24l'(3l' \pm 2) + 5$. Clearly, $p(l') \equiv 5 \pmod{8}$. Since coefficient of $l'^2$ is positive, $p(l')$ is irreducible over $\mathbb{Z}$ and gcd of $p(1)$ and $p(2)$ is one. To proceed further, we state a famous conjecture of Bunyakovsky \cite{BU1857}.
\begin{conj}[Bunyakovsky conjecture]\label{bc}
Let f(x) be a polynomial of one variable with integer coefficient having the following property.
\begin{itemize}
\item[(i)]  The leading coefficient must be positive;
\item[(ii)] Polynomial must be irreducible over the integers;
\item[(iii)] The values $f(1), f(2), f(3), \cdots$ must have no common factor.
\end{itemize}
Then there are infinitely many primes in the sequence $f(1), f(2), f(3), \cdots.$
\end{conj}
Applying this conjecture, we conclude that there exist infinitely many primes of form $24l'(3l' \pm 2) + 5$. Therefore, there are infinitely many $d$'s for which Conjecture \ref{Con1.1} fails.  

\subsection{Toy examples} Here we will discuss a couple of explicit counter-examples to Conjecture \ref{1.1}. 
\subsection*{Example 1} Assume that $d = 10$ and $n = 26 + 6\sqrt{10}$. 
Since $(13)^2 - 10(3)^2 = 79$ is a prime number and $79 \equiv 3 \pmod4$, so that as discussed earlier,  $n$ can not be represented as a difference of two squares in $\mathbb{Z}[\sqrt{10}]$. However, by Theorem \ref{1.1} we can   generate
infinitely many quadruples in $\mathbb{Z}[\sqrt{10}]$ with the property $D(n)$. One of them is 
 $\{19 + 6\sqrt{10}, -8 + 6\sqrt{10}, 35 + 18\sqrt{10},  35 + 42\sqrt{10} \}$ since 
$$
(19 + 6\sqrt{10})(-8 + 6\sqrt{10}) + n = (12 + 3\sqrt{10})^2,
$$ 
$$
(19 + 6\sqrt{10})(35 + 18\sqrt{10}) + n = (31 +  9\sqrt{10})^2, 
$$ 
$$
(19 + 6\sqrt{10})(35 + 42\sqrt{10}) + n = (39 + 13\sqrt{10})^2, 
$$ 
$$
(-8 + 6\sqrt{10})(35 + 18\sqrt{10}) + n = (4 + 9\sqrt{10})^2,
$$
$$
(-8 + 6\sqrt{10})(35 + 42\sqrt{10}) + n = (-4 + 15\sqrt{10})^2,
$$
$$
(35 + 18\sqrt{10})(35 + 42\sqrt{10}) + n = (39 + 27\sqrt{10})^2. 
$$ 
\subsection*{Example 2} Let $d = 58$ and $n = 18 + 2\sqrt{58}$. 
As $(9)^2 - 58(1)^2 = 23$ is a prime number and $23 \equiv 3 \pmod4$, so that $n$ can not be represented as a difference of two squares in $\mathbb{Z}[\sqrt{58}]$. On the other hand, utilizing the case 5 of the proof of  Theorem \ref{1.1} we find a quadruple in $\mathbb{Z}[\sqrt{58}]$ with the property $D(n)$, viz.  
$\{ 19603 + 2574\sqrt{58},  543627 -70094\sqrt{58},  543616 - 70094\sqrt{58},  2154883 - 282950\sqrt{58}\}$ since
$$
(19603 + 2574\sqrt{58})(543627 -70094\sqrt{58}) + n =  (9807 + 1287\sqrt{58} )^2,
$$ 
$$
(19603 + 2574\sqrt{58})(543616 - 70094\sqrt{58}) + n  = ( 9796 + 1287\sqrt{58})^2,
$$ 
$$
(19603 + 2574\sqrt{58})(2154883 - 282950\sqrt{58}) + n = (- 3 + \sqrt{58})^2,
$$ 
$$
(543627 -70094\sqrt{58})(543616 - 70094\sqrt{58}) + n = (- 533820 + 71381\sqrt{58})^2,
$$ 
$$
(543627 -70094\sqrt{58})(2154883 - 282950\sqrt{58}) + n  = ( 1077447 -141475\sqrt{58})^2,
$$ 
$$
(543616 - 70094\sqrt{58})(2154883 - 282950\sqrt{58}) + n  = (- 1077436 + 141475\sqrt{58})^2.
$$ 

\section*{Acknowledgements}
The authors thank the anonymous referee for his/her valuable remarks/suggestions that immensely improved the results as well as the presentation of the paper. The third author would like to appreciate the hospitality provided by Kerala School of Mathematics, Kozhikode, Kerala, where the a part of the work was done. The third author acknowledges SERB MATRICS grant (No. MTR/2021/00762), Govt. of India.

\subsection*{Declaration} The authors declare that there are no conflict of interests.

\subsection*{Data Availability Statement} This manuscript has no associate data.


\begin{thebibliography}{99}
\bibitem{MR2004} F. S. Abu Muriefah and A. Al Rashed, \textit{Some Diophantine quadruples in the ring $\mathbb{Z}[\sqrt{-2}]$}, Math. Commun. \textbf{9} (2004), 1--8.

\bibitem{BD1969} A. Baker and H. Davenport, {\it The equations $3x^2 - 2 = y^2$ and $8x^2 - 7 = z^2$}, Quart. J. Math. Oxford Ser. (2) {\bf 20} (1969), 129--137.


\bibitem{BCM2020} N. C. Bonciocat, M. Cipu and M. Mignotte, {\it There is no Diophantine $D(-1)$-quadruple}, J. London Math. Soc. \textbf{105} (2022), 63-69.


\bibitem{BR1985} E. Brown, {\it Sets in which $xy + k$ is always a square}, Math. Comp. {\bf 45} (1985), 613--620.

\bibitem{BU1857} V. Buniakovsky, {\it Sur les diviseurs numeriques invariables des fonctions rationnelles entieres}, Mem Acad. Sci. St Petersburg {\bf 6} (1857), 305--329.

\bibitem{CGH-2022} K. Chakraborty, S. Gupta and A. Hoque, {\it Diophantine triples with the property $D(n)$ for distinct $n$'s}, Mediterr. J. Math. (to appear). 


\bibitem{DU1993} A. Dujella, {\it Generalization of a problem of Diophantus}, Acta Arith. \textbf{65} (1993), 15--27.

\bibitem{DU1997} A. Dujella, {\it The problem of Diophantus and Davenport for Gaussian integers}, Glas. Mat. Ser. III {\bf 32} (1997), 1--10.

\bibitem{DU2004} A. Dujella, {\it There are only finitely many Diophantine quintuples}, J. Reine Angew. Math. \textbf{566} (2004), 183--214.

\bibitem{DF2007} A. Dujella and Franu\v si\'c, \textit{On differences of two squares in some quadratic fields},  Rocky Mountain J. Math. \textbf{37} (2007), 429--453.

\bibitem{DU21}  A. Dujella, \textit{Number Theory}, \v Skolska knjiga, Zagreb, 2021.

\bibitem{EFF2014} C. Elsholz, A. Filipin and Y. Fujita, {\it On Diophantine quintuples and $D(-1)$-quadruples}, Monats. Math. {\bf 175} (2014), 227--239.

\bibitem{FR2004} Z. Franu\v si\'c, \textit{Diophantine quadruples in the ring $\mathbb{Z}[\sqrt{2}]$}, Math. Commun. \textbf{9} (2004), 141--148.

\bibitem{FR2008} Z. Franu\v si\'c, \textit{Diophantine quadruples in $\mathbb{Z}[\sqrt{4k + 3}]$}, Ramanujan J. \textbf{17} (2008), 77--88.

\bibitem{FR2009} Z. Franu\v si\'c, \textit{A Diophantine problem in $\mathbb{Z}[\sqrt{(1 + d)/2}]$}, Studia Sci. Math. Hungar. \textbf{46} (2009), 103--112.

\bibitem{FR2013} Z. Franu\v si\'c, \textit{Diophantine quadruples in the ring of integers of the pure cubic field $\mathbb{Q}(\sqrt[3]{2})$}, Miskolc Math.
Notes \textbf{14} (2013), 893--903.

\bibitem{FJ2019} Z. Franu\v si\'c and B. Jadrijevi\' c, \textit{$D(n)$-quadruples in the ring of integers of $\mathbb{Q}(\sqrt{2}, \sqrt{3})$}, Math. Slovaca \textbf{69} (2019), 1263--1278.


\bibitem{FS2014} Z. Franu\v si\'c and I. Soldo, \textit{The problem of Diophantus for integers of $\mathbb{Q}(\sqrt{-3})$}, Rad Hrvat. Akad. Znan. Umjet. Mat. Znan. \textbf{18} (2014), 15--25.

\bibitem{HT2011} B. He and A. Togb\'e, {\it On the $D(-1)$-triple $\{1, k^2+1, k^2+2k + 2\}$ and its unique $D(1)$-extension}, J. Number Theory {\bf 131} (2011), 120--137.

\bibitem{HTZ2019} B. He, A. Togb\'e and V. Ziegler, {\it There is no Diophantine quintuple}, Trans. Amer. Math. Soc. \textbf{371} (2019), 6665--6709.

\bibitem{LEM} F. Lemmermeyer, \textit{Higher descent on Pell conics I. From Legendre to Selmer}, preprint. \url{ math.NT/0311309}.

\bibitem{MA2012} Lj. Juki\'c Mati\'c, {\it Non-existence of certain Diophantine quadruples in rings of integers of pure cubic fields}, Proc. Japan Acad. Ser. A  Math. Sci. {\bf 88} (2012), no. 10, 163--167.

\bibitem{ME1895} F. Mertens, {\it \"{U}eber Dirichletsche Reihen}, Ak. Wiss. Wien. (Math.) {\bf 104} (1895), 1093--1153. 

\bibitem{ME1888} A. Meyer,  {\it \"{U}eber einen Satz von Dirichlet}, J. Reine Angew. Math. {\bf 103} (1888), 98--117.


\bibitem{MR1985} S. P. Mohanty and M. S. Ramamsamy, {\it On $P_{r,k}$ sequences}, Fibonacci Quart. {\bf 23} (1985), 36--44.

\bibitem{SO2013} I.  Soldo, \textit{On the existence of Diophantine quadruples in $\mathbb{Z}[\sqrt{-2}]$}, Miskolc Math. Notes \textbf{14} (2013), 265--277.

\end{thebibliography}
\end{document}